\documentclass[12pt,draft]{amsart}
\usepackage{amssymb}
\usepackage{enumitem}
\usepackage{graphicx}
\def\pmod #1{\ ({\rm{mod}}\ #1)}
\def\Z{\Bbb Z}

\def\bg{\bigg}
\def\({\bg(}
\def\){\bg)}

\def\sgn{{\rm sgn}}

\def\sgn{{\rm sgn}}

\def\Ack{\medskip\noindent {\bf Acknowledgments}}
\theoremstyle{plain}
\newtheorem{theorem}{Theorem}

\newtheorem{lemma}{Lemma}
\newtheorem{corollary}{Corollary}

\theoremstyle{definition}

\theoremstyle{remark}
\newtheorem{remark}{Remark}

\makeatletter
\@namedef{subjclassname@2010}{%
  \textup{2010} Mathematics Subject Classification}
\makeatother
 \vspace{4mm}

\begin{document}
 \baselineskip=17pt
\hbox{} {}
\medskip
\title[H.-L. Wu and L.-Y. Wang]
{Applications of circulant matrices to determinants involving $k$-th power residues}
\date{}
\author[H.-L. Wu and L.-Y. Wang]{Hai-Liang Wu and Li-Yuan Wang}

\thanks{2020 {\it Mathematics Subject Classification}.
Primary 11C20; Secondary 11L05, 11R29.
\newline\indent {\it Keywords}. determinants, the Legendre symbol, circulant matrices.
\newline \indent  This research was supported by the National Natural Science Foundation of China (Grant No. 12101321 and Grant No. 11971222).}

\address {(Hai-Liang Wu) School of Science, Nanjing University of Posts and Telecommunications, Nanjing 210023, People's Republic of China}
\email{\tt whl.math@smail.nju.edu.cn}

\address {(Li-Yuan Wang) School of Physical and Mathematical Sciences, Nanjing Tech University, Nanjing 211816, People's Republic of China}
\email{\tt wly@smail.nju.edu.cn}

\begin{abstract}
In this paper, by the tools of circulant matrices and hyperelliptic curves over finite fields, we study some arithmetic properties of certain determinants involving the Legendre symbols and $k$-th residues.
\end{abstract}
\maketitle
\section{Introduction}
Let $n$ be an arbitrary positive integer and let $R$ be a commutative ring. For every $n\times n$ matrix $M=[a_{ij}]_{1\le i,j\le n}$ with $a_{ij}\in R$, we use the symbol $\det M$ or $|M|$ to denote the determinant of $M$. Also, given any elements $b_0,b_1,\cdots,b_{n-1}\in R$, the {\it circulant matrix} of the $n$-tuple $(b_0,\cdots,b_{n-1})$ is defined by an $n\times n$ matrix over $R$ whose $(i,j)$-entry is $b_{i-j}$, where the indices are cyclic modulo $n$. We also denote this matrix by $C(b_0,b_1,\cdots,b_{n-1})$. Readers may refer to the survey paper \cite{Kra} for more results on circulant matrices.
\subsection{Circulant Matrices involving The Legendre Symbols}
Circulant matrices have many applications in both number theory and combinatorics. Let $p$ be an odd prime and let $\chi(\cdot)$ be a multiplicative character modulo $p$. Carlitz \cite{carlitz} first investigated the following circulant matrix:
$$
C(c_0,c_1,\cdots,c_{p-1}):=\bigg[\mu+\chi(i-j)\bigg]_{1\le i,j\le p-1}\ \ \ \ \ (\mu\in\mathbb{C}),
$$
where $c_i=\mu+\chi(i)$ for $0\le i\le p-1$. Carlitz \cite[Theorem 4]{carlitz} determined the characteristic polynomial of this circulant matrix. In particular, when $\chi(\cdot)=(\frac{\cdot}{p})$ is the Legendre symbol, the characteristic polynomial of the matrix $[\mu+(\frac{i-j}{p})]_{1\le i,j\le p-1}$ is
$$F_{\mu}(t)=(t^2-(-1)^{(p-1)/2}p)^{(p-3)/2}(t^2-(p-1)\mu-(-1)^{(p-1)/2}).$$
Later Chapman \cite{chapman,evil} and Vsemirnov \cite{M1,M2} also studied many variants of Carlitz's results involving the Legendre symbols.

Let $p=2n+1$ be an odd prime. Recently, Sun \cite{ffadeterminant} studied the following determinant:
$$
S(d,p):=\det\bigg[\left(\frac{i^2+dj^2}{p}\right)\bigg]_{1\le i,j\le n},
$$
where $d\in\mathbb{Z}$ with $p\nmid d$. Sun \cite[Theorem 1.2(iii) and Theorem 1.3(i)]{ffadeterminant} proved that $-S(d,p)$ is a quadratic residue modulo $p$ whenever $(\frac{d}{p})=1$. For the recent progress on this topic, readers may refer to \cite{Dimitry,Wu}. On the other hand, Sun also investigated some global properties of this determinant. Sun conjectured that $-S(1,p)$ is an integral square if $p\equiv 3\pmod 4$. Later, by using a sophisticated matrix decomposition, Alekseyev and Krachun proved this conjecture. Also, in the case $p\equiv1\pmod 4$, writing $p=a^2+4b^2$ with $a,b\in\Z$ and $a\equiv 1\pmod 4$, Cohen, Sun and Vsemirnov (see \cite[Remark 4.2]{ffadeterminant}) conjectured that $S(1,p)/a$ is an integral square. This conjecture was later proved by the first author \cite[Theorem 3]{Wu1}.

It is worth to state here that $S(d,p)$ is indeed a determinant of certain circulant matrix. In fact, fix a primitive root $g$ modulo $p$. Then it is clear that $S(d,p)$ is equal to
\begin{align*}
\det\bigg[\left(\frac{g^{2i}+dg^{2j}}{p}\right)\bigg]_{0\le i,j\le n-1}
&=\det\bigg[\left(\frac{g^{2(i-j)}+d}{p}\right)\bigg]_{0\le i,j\le n-1}\\
&=\det C(s_0,s_1,\cdots,s_{n-1}),
\end{align*}
where $s_i=(\frac{g^{2i}+d}{p})$ for $0\le i\le n-1$.

Motivated by Sun's determinant $S(d,p)$, in this paper we study some determinants concerning $k$-th power residues. Let $p$ be an odd prime and let $k\ge 2$ be an integer dividing $p-1$. Write $p=km+1$ and let
$$
0<\alpha_1<\alpha_2<\cdots<\alpha_m<p
$$
be all the $k$-th power residues modulo $p$ in the interval $(0,p)$. We consider the following matrix:
\begin{equation}\label{Eq. Definition of the matrix Wp(k)}
W_p(k):=\bigg[\left(\frac{\alpha_i+\alpha_j}{p}\right)\bigg]_{1\le i,j\le m}.
\end{equation}
To state our results, we first introduce some notations. Let $\mathbb{F}_p$ denote the finite field of $p$ elements. Let $\mathcal{C}_{p,k,\psi}$ and $\mathcal{C}_{p,k,\phi}$ be the curves over $\mathbb{F}_p$ defined by the equations $y^2=x^k+1$ and $y^2=x(x^k+1)$ respectively. Also, we define $a_p(k)$ and $b_p(k)$ by
\begin{equation}\label{Eq. Definition of ap(k)}
p+1-a_p(k)=\#\{(x,y)\in\mathbb{F}_p\times\mathbb{F}_p:\ y^2=x^k+1\}\cup\{\infty\},
\end{equation}
and
\begin{equation}\label{Eq. Defintition of bp(k)}
p+1-b_p(k)=\#\{(x,y)\in\mathbb{F}_p\times\mathbb{F}_p:\ y^2=x(x^k+1)\}\cup\{\infty\},
\end{equation}
where $\#S$ denotes the cardinality of a set $S$.

When $k$ is even, as a generalization of Sun's determinant $S(1,p)$, by using the theory of circulant matrices we can obtain the following result:
\begin{theorem}\label{Thm. A}
Let $p$ be an odd prime and let $k\ge 2$ be an even integer dividing $p-1$. Then the following results hold.

{\rm (i)} If $m$ is odd, then
$$
\det W_p(k)=-(a_p(k)+1)u_p(k)^2/k
$$
for some $u_p(k)\in\mathbb{Z}$.

{\rm (ii)} If $m$ is even, then
$$
\det W_p(k)=(a_p(k)+1)b_p(k)v_p(k)^2/k^2
$$
for some $v_p(k)\in\mathbb{Z}$.
\end{theorem}
\begin{remark}\label{Rem. of Thm. A}
{\rm (i)} When $k=2$ and $p\equiv3\pmod4$, it is easy to see that $a_p(2)=1$. This implies that $-\det W_p(2)=-S(1,p)$ is an integral square, which also confirms the conjecture of Sun.

{\rm (2)} When $k=2$ and $p\equiv 1\pmod4$ with $p=a^2+4b^2$, where $a\equiv1\pmod4$, it is known that $a_p(2)=1$ and $b_p(2)=2a$ (cf. \cite[Theorem 6.2.9]{Berndt}). This gives that $\det W_p(2)/a=S(1,p)/a$ is an integral square, which coincides with the result in \cite[Theorem 3]{Wu1}.
\end{remark}
Now we consider the case that $k$ is odd.

To state our next result, we need to introduce some notations. Fix a primitive root $g$ modulo $p$. Let $E_{p,k,1}$ and $E_{p,k,g}$ be the hyperelliptic curves over $\mathbb{F}_p$ defined by the equations $y^2=x(x^{2k}+1)$ and $y^2=x(x^{2k}+g^k)$ respectively. We also define $c_p(k)$ and $d_p(k)$ by
\begin{equation}\label{Eq. Definition of cp(k)}
p+1-c_p(k):=\#\{(x,y)\in\mathbb{F}_p\times\mathbb{F}_p:\ y^2=x(x^{2k}+1)\}\cup\{\infty\},
\end{equation}
and
\begin{equation}\label{Eq. Definition of dp(k)}
p+1-d_p(k):=\#\{(x,y)\in\mathbb{F}_p\times\mathbb{F}_p:\ y^2=x(x^{2k}+g^k)\}\cup\{\infty\}.
\end{equation}

Now we state our next result.

\begin{theorem}\label{Thm. B}
{\rm (i)} Let $p\equiv 1\pmod 4$ be a prime and let $k\ge 2$ be an odd integer dividing $p-1$. Then we have
$$
\det W_p(k)=\frac{z_p(k)^2}{4k^2}\left(c_p(k)^2+d_p(k)^2\right)
$$
for some $z_p(k)\in\mathbb{Z}$.

{\rm (ii)} Let $p\equiv 3\pmod 4$ be a prime and let $k\ge 2$ be an odd integer dividing $p-1$. Then $-\det W_p(k)$ is an integral square.
\end{theorem}
When $k=3$ we have the following result:
\begin{corollary}\label{Coro1. of Thm. B}
Suppose that $p\equiv1\pmod{12}$ is a prime and write $p=c^2+9d^2$ with $c,d\in\mathbb{Z}$. Then

{\rm (i)} $\det W_p(3)/(c^2+d^2)$ is an integral square.

{\rm (ii)} Moreover, if $p\nmid\det W_p(3)$, then
$$
\left(\frac{\det W_p(3)}{p}\right)=\left(\frac{2}{p}\right).
$$
\end{corollary}
\begin{remark}\label{Rem. of corollary}
By our computations, there are primes $p\equiv1\pmod {12}$ such that $p\mid\det W_p(3)$. For example, $1117,1129,1381,1597,1861,2557,2749$ are all primes $p\equiv1\pmod{12}$ less than $3000$ such that $p\mid W_p(3)$.
\end{remark}
We also have the following corollary:
\begin{corollary}\label{Coro2. of Thm. B}
{\rm (i)} Let $p\equiv 1\pmod 4$ be a prime and let $k\ge 2$ be an odd integer dividing $p-1$. Then $\det W_p(k)\ge0$.

{\rm (ii)} Let $p\equiv 3\pmod 4$ be a prime and let $k\ge 2$ be an odd integer dividing $p-1$. Then $\det W_p(k)\le 0$.
\end{corollary}

\subsection{Determinants of the form $\det[\frac{1}{\alpha_i+\alpha_j}]_{1\le i,j\le m}$}
Let $p$ be an odd prime. In 2019, Sun \cite{ffadeterminant} also studied the following determinant:
$$
A_p:=\det\bigg[\frac{1}{i^2+j^2}\bigg]_{1\le i,j\le (p-1)/2}.
$$
When $p\equiv3\pmod4$, Sun \cite[Theorem 1.4(ii)]{ffadeterminant} showed that
$$
A_p\equiv \left(\frac{2}{p}\right)\pmod p.
$$
In \cite[Remark 1.3]{ffadeterminant} Sun also conjectured that if $p\equiv 2\pmod 3$ is odd, then $2B_p$ is a quadratic residue modulo $p$, where
$$
B_p:=\det\bigg[\frac{1}{i^2-ij+j^2}\bigg]_{1\le i,j\le p-1}.
$$
This conjecture was later confirmed in \cite{WSN}. Let the notations be as the above. Inspired by the above work, we consider the matrix:
\begin{equation}\label{Eq. Definition of the matrix Ip(k)}
I_p(k):=\bigg[\frac{1}{\alpha_i+\alpha_j}\bigg]_{1\le i,j\le m}.
\end{equation}
As a generalization of Sun's determinant $\det A_p$, we obtain the following result:
\begin{theorem}\label{Thm. C}
Let $p$ be an odd prime and let $k\ge2$ be an even integer dividing $p-1$. Write $p=km+1$. Suppose that $-1$ is not a $k$-th power residue modulo $p$. Then we have
$$
\det I_p(k)\equiv \frac{(-1)^{\frac{m+1}{2}}}{(2k)^m}\pmod p.
$$
\end{theorem}
\begin{remark}\label{Rem. of Thm. C}
When $p\equiv3\pmod4$ and $k=2$, by the above theorem
$$
\det I_p(2)\equiv (-1)^{\frac{p+1}{4}}=\left(\frac{2}{p}\right)\pmod p.
$$
This coincides with Sun's result \cite[Theroem 1.4(ii)]{ffadeterminant}.
\end{remark}
The outline of this paper is as follows. We will prove Theorems \ref{Thm. A}--\ref{Thm. B} and their corollaries in Section 2. The proof of Theorem \ref{Thm. C} will be given in section 3.

\section{Proofs of Theorems \ref{Thm. A}--\ref{Thm. B}}
Recall that $C(a_0,\cdots,a_{n-1})$ denotes the circulant matrix of the $n$-tuple $(a_0,\cdots,a_{n-1})$. We begin with the following lemma which is the key element of our proofs (cf. \cite[Lemma 3.4]{Wu}).
\begin{lemma}\label{Lem. The key lemma}
Let $R$ be a commutative ring. Let $n$ be a positive integer. Let $a_0,a_1,\cdots,a_{n-1}\in R$ such that
\begin{equation}\label{Eq. condition for ai}
a_i=a_{n-i}\ \ \ \ \ \text{for each}\ 1\le i\le n-1.
\end{equation}
If $n$ is even, then there exists an element $u\in R$ such that
\begin{equation}\label{Eq. result of the key lemma}
\det C(a_0,a_1,\cdots,a_{n-1})=\(\sum_{i=0}^{n-1}a_i\)\(\sum_{i=0}^{n-1}(-1)^ia_i\)u^2.
\end{equation}
If $n$ is odd, then there exists an element $v\in R$ such that
\begin{equation}\label{Eq. result2 of the key lemma}
\det C(a_0,a_1,\cdots,a_{n-1})=\(\sum_{i=0}^{n-1}a_i\)v^2.
\end{equation}
\end{lemma}

{\bf Proof of Theorem \ref{Thm. A}.} Fix a primitive root $g$ modulo $p$. As $k$ is even, we have
\begin{align*}
\det W_p(k)=\det\bigg[\left(\frac{1+\alpha_i/\alpha_j}{p}\right)\bigg]_{1\le i,j\le m}
&=\det\bigg[\left(\frac{1+g^{k(i-j)}}{p}\right)\bigg]_{0\le i,j\le m-1}\\
&=\det C(e_0,e_1,\cdots,e_{m-1}),
\end{align*}
where $e_i=(\frac{1+g^{ki}}{p})$ for $0\le i\le m-1$. Clearly $e_0,\cdots,e_{m-1}$ satisfy the condition (\ref{Eq. condition for ai}). Moreover,
\begin{equation}\label{Eq. a1 in the proof of Thm A}
\sum_{i=0}^{m-1}e_i
=\frac{1}{k}\sum_{x=1}^{p-1}\left(\frac{1+x^k}{p}\right)
=\frac{1}{k}\left(-1+\sum_{x=0}^{p-1}\left(\frac{1+x^k}{p}\right)\right)
=-\left(1+a_p(k)\right)/k
\end{equation}
where $a_p(k)$ is defined by (\ref{Eq. Definition of ap(k)}). Also,
\begin{equation}\label{Eq. a2 in the proof of Thm A}
\sum_{i=0}^{m-2}(-1)^ie_i=\frac{1}{k}\sum_{x=1}^{p-1}\left(\frac{x^k+1}{p}\right)\left(\frac{x}{p}\right)
=-b_p(k)/k,
\end{equation}
where $b_p(k)$ is defined by (\ref{Eq. Defintition of bp(k)}). Combining Lemma \ref{Lem. The key lemma} with (\ref{Eq. a1 in the proof of Thm A}) and (\ref{Eq. a2 in the proof of Thm A}), one can get the desired result.
\qed

Now we turn to the proof of Theorem \ref{Thm. B}. We first need the following known result in linear algebra.
\begin{lemma}\label{Lemma eigenvalues}
Let $M$ be an $n\times n$ complex matrix. Let $\lambda_1,\cdots,\lambda_n$ be complex numbers and let ${\bf u}_1,\cdots, {\bf u}_n$ be $m$-dimensional column vectors. Suppose that $M{\bf u}_i=\lambda_i{\bf u}_i$ for each $1\le i\le n$ and that ${\bf u}_1,\cdots, {\bf u}_n$ are linearly independent. Then $\lambda_1,\cdots,\lambda_n$ are exactly all the eigenvalues of $M$ (counting multiplicities).
\end{lemma}

To state our proof, we introduce the following notations. Let $\widehat{\mathbb{F}_p^{\times}}$ denote the cyclic group of all multiplicative characters of $\mathbb{F}_p$ and let $\chi_p(\cdot)$ be a generator of $\widehat{\mathbb{F}_p^{\times}}$. For any matrix $M$, we use the symbol $M^T$ to denote the transpose of $M$.

Now we are in a position to proof our main result.

{\bf Proof of Theorem \ref{Thm. B}.} Recall that $k\ge2$ is an odd integer dividing $p-1$ and $p=km+1$.

(i) We first consider the case $p\equiv1\pmod4$. Clearly the elements $\alpha_1\ {\rm mod}\ p ,\cdots, \alpha_m\ {\rm mod}\ p$ are exactly $m$ distinct roots of the polynomial $X^m-1$ over $\mathbb{F}_p=\mathbb{Z}/p\mathbb{Z}$. We therefore have
\begin{equation}\label{Eq. b1 in the proof of Thm B}
X^m-1\equiv \prod_{i=1}^{m}\left(X-\alpha_i\right)\pmod p.
\end{equation}
By (\ref{Eq. b1 in the proof of Thm B}) we obtain
\begin{equation}\label{Eq. b2 in the proof of Thm B}
\prod_{j=1}^m\alpha_j\equiv -1^{m+1}=-1\pmod p.
\end{equation}
By (\ref{Eq. b2 in the proof of Thm B}) it is easy to see that $\det W_p(k)$ is equal to
\begin{equation}\label{Eq. b3 in the proof of Thm B}
\left(\frac{-1}{p}\right)
\det\bigg[\left(\frac{\alpha_i+\alpha_j}{p}\right)\bigg]_{1\le i,j\le m}
=\det\bigg[\left(\frac{\alpha_i+\alpha_j}{p}\right)\left(\frac{\alpha_j}{p}\right)\bigg]_{1\le i,j\le m}.
\end{equation}
Now we first determine all the eigenvalues of the matrix:
$$
W^*_p(k):=\bigg[\left(\frac{\alpha_i+\alpha_j}{p}\right)\left(\frac{\alpha_j}{p}\right)\bigg]_{1\le i,j\le m}.
$$
For each $1\le r\le m$, we have
\begin{align*}
\sum_{j=1}^{m}\left(\frac{\alpha_i+\alpha_j}{p}\right)\left(\frac{\alpha_j}{p}\right)\chi_p^r(\alpha_j)
&=\sum_{j=1}^{m}\left(\frac{1+\alpha_j/\alpha_i}{p}\right)\left(\frac{\alpha_j/\alpha_i}{p}\right)
\chi_p(\alpha_j/\alpha_i)\chi_p^r(\alpha_i)\\
&=\sum_{j=1}^{m}\left(\frac{1+\alpha_j}{p}\right)\left(\frac{\alpha_j}{p}\right)
\chi_p^r(\alpha_j)\chi_p^r(\alpha_i).
\end{align*}
This implies that for each $1\le r\le m$ we have
$$
W^*_p(k){\bf v}_r=\lambda_r{\bf v}_r,
$$
where
$$
\lambda_r=\sum_{j=1}^{m}\left(\frac{1+\alpha_j}{p}\right)\left(\frac{\alpha_j}{p}\right)
\chi_p^r(\alpha_j),
$$
and
$$
{\bf v}_r=\left(\chi_p^r(\alpha_1),\cdots,\chi_p^r(\alpha_m)\right)^T.
$$
Note that
$$
\left| \begin{array}{cccccccc}
\chi_p^1(\alpha_1) & \chi_p^2(\alpha_1) & \ldots  & \chi_p^m(\alpha_1)  \\
\chi_p^1(\alpha_2) & \chi_p^2(\alpha_2) &  \ldots  &  \chi_p^n(\alpha_2)\\
\vdots & \vdots & \ddots  & \vdots    \\
\chi_p^1(\alpha_n)& \chi_p^2(\alpha_n) & \ldots &  \chi_p^m(\alpha_m)\\
\end{array} \right|=\pm\prod_{1\le i<j\le m}\(\chi_p(\alpha_j)-\chi_p(\alpha_i)\)\ne0.
$$
Hence the vectors ${\bf v}_1,\cdots, {\bf v}_m$ are linearly independent. Now by Lemma \ref{Lemma eigenvalues} the numbers $\lambda_1,\cdots,\lambda_m$ are exactly all the eigenvalues of $W^*_p(k)$ (counting multiplicities).

When $r=m$,
\begin{equation}\label{Eq. b4 in the proof of Thm B}
\lambda_m=\sum_{j=1}^m\left(\frac{1+\alpha_j}{p}\right)\left(\frac{\alpha_j}{p}\right)
=\frac{1}{k}\sum_{x=1}^{p-1}\left(\frac{1+x^k}{p}\right)\left(\frac{x}{p}\right).
\end{equation}
When $r=m/2$,
\begin{equation}\label{Eq. b4 in the proof of Thm B}
\lambda_{m/2}=\sum_{j=1}^m\left(\frac{1+\alpha_j}{p}\right)
=\frac{1}{k}\sum_{x=1}^{p-1}\left(\frac{1+x^k}{p}\right).
\end{equation}
By \cite[Proposition 6.1.7]{Berndt} we have
\begin{equation}\label{Eq. b5 in the proof of Thm B}
\lambda_m=\lambda_{m/2}.
\end{equation}
In addition, when $1\le r\le m/2-1$ it is clear that $\overline{\lambda_r}=\lambda_{m-r}$, where $\bar{z}$ denotes the complex conjugation of a complex number $z$. Combining this with (\ref{Eq. b5 in the proof of Thm B}), we have
\begin{equation}\label{Eq. det Wp(k) is greater than 0}
\det W_p(k)=\det W^*_p(k)=\prod_{r=1}^m\lambda_r
=\lambda_m^2\prod_{1\le r\le m/2-1}\lambda_r\overline{\lambda_r}\ge0.
\end{equation}
Let ${\bf i}\in\mathbb{C}$ be a $4$-th primitive root of unity. Fix a primitive root $g$ modulo $p$. Then
\begin{align*}
\det W^*_p(k)
&=\det\bigg[\left(\frac{\alpha_i+\alpha_j}{p}\right)\left(\frac{\alpha_j}{p}\right){\bf i}^{i-j}\bigg]_{1\le i,j\le m}\\
&=\det\bigg[\left(\frac{1+g^{k(i-j)}}{p}\right){\bf i}^{i-j}\bigg]_{0\le i,j\le m-1}\\
&=\det C(\omega_0,\cdots,\omega_{m-1}),
\end{align*}
where $\omega_r=(\frac{1+g^{kr}}{p}){\bf i}^r$ for $0\le r\le m-1$. One can verify that  $\omega_0,\cdots,\omega_{m-1}$ satisfy the condition (\ref{Eq. condition for ai}). Fix a multiplicative character $\psi\in\widehat{\mathbb{F}_p^{\times}}$ of order $4$ with $\psi(g)={\bf i}$. Then
\begin{equation*}
\sum_{r=0}^{m-1}\omega_r
=\sum_{r=0}^{m-1}\left(\frac{1+g^{kr}}{p}\right)\psi(g^r)
=\frac{1}{k}\sum_{r=0}^{p-2}\left(\frac{1+g^{kr}}{p}\right)\psi(g^r).
\end{equation*}
One can also verify the following equalities:
\begin{align*}
\sum_{r=0}^{p-2}\left(\frac{1+g^{kr}}{p}\right)\psi(g^r)
&=\sum_{r=0}^{\frac{p-3}{2}}\left(\frac{1+g^{2kr}}{p}\right)\left(\frac{g^r}{p}\right)
+{\bf i}\sum_{r=0}^{\frac{p-3}{2}}\left(\frac{1+g^{2kr}g^k}{p}\right)\left(\frac{g^r}{p}\right)\\
&=\frac{1}{2}\sum_{x=1}^{p-1}\left(\frac{1+x^{2k}}{p}\right)\left(\frac{x}{p}\right)
+\frac{1}{2}{\bf i}\sum_{x=1}^{p-1}\left(\frac{1+x^{2k}g^k}{p}\right)\left(\frac{x}{p}\right)\\
&=\frac{1}{2}\sum_{x=1}^{p-1}\left(\frac{1+x^{2k}}{p}\right)\left(\frac{x}{p}\right)
+\frac{1}{2}{\bf i}\sum_{x=1}^{p-1}\left(\frac{g^k+x^{2k}}{p}\right)\left(\frac{x}{p}\right)\\
&=-\left(c_p(k)+{\bf i}d_p(k)\right)/2,
\end{align*}
where $c_p(k)$ and $d_p(k)$ are defined by (\ref{Eq. Definition of cp(k)}) and (\ref{Eq. Definition of dp(k)}) respectively. Hence
\begin{equation}\label{Eq. sum of omega r in the proof of Thm. B}
\sum_{r=0}^{m-1}\omega_r=\frac{-1}{2k}\left(c_p(k)+{\bf i}d_p(k)\right).
\end{equation}
With essentially the same method, one can also verify that
\begin{equation}\label{Eq. sum of +- omega r in the proof of Thm. B}
\sum_{r=0}^{m-1}(-1)^r\omega_r=\frac{-1}{2k}\left(c_p(k)-{\bf i}d_p(k)\right).
\end{equation}
If $\det W_p(k)=0$, then one can get the desired result directly. Suppose now that $\det W_p(k)\ne0$. By (\ref{Eq. det Wp(k) is greater than 0}) we have $\det W_p(k)>0$ under this assumption.
Combining Lemma \ref{Lem. The key lemma} with (\ref{Eq. sum of omega r in the proof of Thm. B}) and (\ref{Eq. sum of +- omega r in the proof of Thm. B}), there exists an element $z_p(k)\in\mathbb{Z}[{\bf i}]$ such that
$$
\det W_p(k)=\det W_p^*(k)=\frac{z_p(k)^2}{4k^2}\left(c_p(k)^2+d_p(k)^2\right).
$$
As $\det W_p(k)\in\mathbb{Z}$ and $\det W_p(k)>0$, the number $z_p(k)$ must be an integer. This completes the proof of (i).

(ii) We now consider the case $p\equiv3\pmod4$. As $k$ is odd, it is clear that
$$
-\alpha_1\ {\rm mod}\ p,\cdots, -\alpha_{m}\ {\rm mod}\ p
$$
is a permutation $\pi$ of the sequence
$$
\alpha_1\ {\rm mod}\ p,\cdots, \alpha_{m}\ {\rm mod}\ p,
$$
and clearly
$$
\sgn(\pi)\equiv\prod_{1\le i<j\le m}\frac{-\alpha_j-(-\alpha_i)}{\alpha_j-\alpha_i}=(-1)^{m(m-1)/2}\pmod p,
$$
where $\sgn(\pi)$ is the sign of $\pi$.
When $p\equiv 3\pmod 4$ and $k$ is odd, since $m\equiv 2\pmod 4$, the number $\det W_p(k)$ is equal to
$$
\sgn(\pi)\det\bigg[\left(\frac{\alpha_i-\alpha_j}{p}\right)\bigg]_{1\le i,j\le m}
=-\det\bigg[\left(\frac{\alpha_i-\alpha_j}{p}\right)\bigg]_{1\le i,j\le m}.
$$
Clearly the matrix $M_p:=[(\frac{\alpha_i-\alpha_j}{p})]_{1\le i,j\le m}$ is a skew-symmetric matrix, i.e., $M_p^{T}=-M_p$. It is known that the determinant of a skew-symmetric matrix of even order with integer entries is always an integral square (cf. \cite[Proposition 2.2]{Stem}). This implies that $-\det W_p(k)$ is an integral square.

In view of the above, we have completed the proof.
\qed

{\bf Proof of Corollary \ref{Coro1. of Thm. B}.} (i) Let $k=3$ and $p\equiv1\pmod{12}$. Write $p=\alpha^2+\beta^2$ with $\alpha,\beta\in\mathbb{Z}$ and $\alpha\equiv-(\frac{2}{p})\pmod 4$. It is known that (cf. \cite[Theorem 6.2.5]{Berndt})
$$
c_p(3)^2=\begin{cases}36\alpha^2&\mbox{if}\ 3\nmid \alpha,\\ 4\alpha^2&\mbox{if}\
3\mid\alpha,\end{cases},\
d_p(3)^2=\begin{cases}4\beta^2&\mbox{if}\ 3\nmid \alpha,\\ 36\beta^2&\mbox{if}\
3\mid\alpha.\end{cases}
$$
Hence if we write $p=c^2+9d^2$ with $c,d\in\Z$, then one can easily verify that
$$
\frac{c_p(3)^2+d_p(3)^2}{36}=c^2+d^2.
$$
By Theorem \ref{Thm. B} we obtain that $\det W_p(3)/(c^2+d^2)$ is an integral square if $p\equiv 1\pmod {12}$.

(ii) If $p\nmid \det W_p(3)$, then
$$
\left(\frac{\det W_p(3)}{p}\right)=\left(\frac{c^2+d^2}{p}\right)
=\left(\frac{8c^2+p}{p}\right)=\left(\frac{2}{p}\right).
$$
This completes the proof.\qed

\section{Proof of Theorem \ref{Thm. C}.}
Recall that
$$
I_p(k)=\bigg[\frac{1}{\alpha_i+\alpha_j}\bigg]_{1\le i,j\le m}.
$$
As $-1$ is not a $k$-th power residue modulo $p$, clearly we have $2\nmid m$. Now we prove our theorem.

{\bf Proof of Theorem \ref{Thm. C}.} By \cite[Theorem 12(5.5)]{K2} we have
\begin{equation}\label{Eq. c1 in the proof of Thm. C}
\det I_p(k)=\frac{\prod_{1\le i<j\le m}\left(\alpha_i-\alpha_j\right)^2}
{\prod_{1\le i\le m}\prod_{1\le j\le m}\left(\alpha_i+\alpha_j\right)}.
\end{equation}
We first consider the numerator. One can verify the following equalities:
\begin{align*}
N_p:=\prod_{1\le i<j\le m}\left(\alpha_i-\alpha_j\right)^2
&=(-1)^{\frac{m(m-1)}{2}}\prod_{1\le i\ne j\le m}\left(\alpha_i-\alpha_j\right)\\
&=(-1)^{\frac{(m-1)}{2}}\prod_{1\le j\le m}\prod_{i\ne j}\left(\alpha_j-\alpha_i\right)\\
&=(-1)^{\frac{(m-1)}{2}}\prod_{1\le j\le m}G'(\alpha_j),
\end{align*}
where $G'(X)$ is the derivative of $G(X)=\prod_{1\le i\le m}(X-\alpha_i)$.
Observe that
\begin{equation}\label{Eq. polynomial congruence in the proof of Thm. C}
G(X)\equiv X^m-1\pmod p.
\end{equation}
Hence $G'(X)\equiv mX^{m-1}\pmod p$ and $\prod_{1\le i\le m}\alpha_i\equiv (-1)^{m+1}=1\pmod p$.

By this we obtain
\begin{align*}
\prod_{1\le i<j\le m}\left(\alpha_i-\alpha_j\right)^2
&=(-1)^{\frac{(m-1)}{2}}\prod_{1\le j\le m}G'(\alpha_j)\\
&\equiv (-1)^{\frac{(m-1)}{2}}m^m\prod_{1\le j\le m}\alpha_j^{m-1}\equiv (-1)^{\frac{(m-1)}{2}}m^m\pmod p.
\end{align*}
Hence
\begin{equation}\label{Eq. the numerator in the proof of Thm. C}
N_p\equiv (-1)^{\frac{(m-1)}{2}}m^m\pmod p.
\end{equation}

Now we turn to the denominator. One can verify the following equalities:
\begin{align*}
D_p:=\prod_{i=1}^m\prod_{j=1}^m\left(\alpha_i+\alpha_j\right)
=\prod_{i=1}^m\alpha_i^m\prod_{j=1}^m\left(1+\alpha_j/\alpha_i\right)
&\equiv \prod_{i=1}^m\prod_{j=1}^m\left(1+\alpha_j\right)\\
&=\prod_{j=1}^m\left(1+\alpha_j\right)^m\pmod p.
\end{align*}
Hence by (\ref{Eq. polynomial congruence in the proof of Thm. C})
\begin{equation}\label{Eq. the denominator in the proof of Thm. C}
D_p\equiv (-1)^mG(-1)^m\equiv 2^m\pmod p.
\end{equation}
Combining (\ref{Eq. the numerator in the proof of Thm. C}) with (\ref{Eq. the denominator in the proof of Thm. C}), we finally obtain
$$
\det I_p(k)\equiv \frac{(-1)^{\frac{m-1}{2}}m^m}{2^m}\equiv\frac{(-1)^{\frac{m+1}{2}}}{(2k)^m}\pmod p.
$$
This completes the proof.
\qed

\Ack\ We would like to thank Prof. Hao Pan for his steadfast encouragement.

\end{document}